\documentclass{article}
\usepackage{amsmath,amsfonts,amssymb, amsthm}
\usepackage[utf8]{inputenc}
\usepackage[T1]{fontenc}
\usepackage[english]{babel}
\usepackage{vmargin}
\usepackage{graphicx}
\newcommand{\dt}{{\delta \! t}}

\newcommand{\IGS}{\includegraphics[height=0.25\hsize,width=0.40\hsize]}
\newcommand{\xR}{\mathbb{R}}
\newcommand{\xN}{\mathbb{N}}
\def\xdif{\,{\rm d}}
\theoremstyle{plain}
\newtheorem{thrm}{\thrmname}[section]
\newtheorem{lmm}[thrm]{\lmname}

\providecommand{\thrmname}{Theorem}
\providecommand{\lmname}{Lemma}
\providecommand{\crllrname}{Corollary}
\providecommand{\prpstnname}{Proposition}
\providecommand{\crtrnname}{Criterion}
\providecommand{\lgrthmname}{Algorithm}
\theoremstyle{definition}

\newtheorem{dfntn}[thrm]{\dfntnname}

\newtheorem{rmrk}[thrm]{\rmrkname}

\providecommand{\dfntnname}{Definition}
\providecommand{\cnjctrname}{Conjecture}
\providecommand{\xmplname}{Example}
\providecommand{\prblmname}{Problem}
\providecommand{\rmrkname}{Remark}
\providecommand{\ntname}{Note}
\providecommand{\clmname}{Claim}
\providecommand{\smmrname}{Summary}
\providecommand{\csnname}{Case}
\providecommand{\bsrvtnname}{Observation}

\begin{document}
\title{Numerical methods for piecewise deterministic Markov processes with boundary}
\author{Ludovic Gouden\`ege
\footnote{Laboratoire d'Analyse et de Math\'ematiques Appliqu\'ees (CNRS - UMR 8050) 
Universit\'e Paris-Est, 5 boulevard Descartes
Champs sur Marne, 77454 Marne-la-Vall\'ee Cedex 2, France}}
\date{September, 2014}
%
%
%
%
\maketitle
\begin{abstract} In this paper is described the general aspect of a numerical method for piecewise deterministic Markov processes with boundary. Under very natural hypotheses, a crucial result about uniqueness of solution of a generalized Kolmogorov equation with respect to a test function space is proved. Next we prove the existence and uniqueness of a positive solution to the finite volume scheme without result about convergence. Finally different models of transmission control protocol window-size processes are simulated to illustrate the efficiency of the numerical method for describing the evolution of the density of a piecewise deterministic Markov process with boundary.
Obviously some technical aspects have been skipped for reader convenience but the full theory will be exposed in a forthcoming paper in collaboration with C. Cocozza-Thivent*, R. Eymard* and M. Roussignol*.
\end{abstract}
\section{Introduction}

Piecewise Deterministic Markov Processes (PDMP) appear in many areas, such as engineering, operations research, biology, economics... One can find the definition and many properties of these processes in the founding book of M.H.A. Davis \cite{DA}. Relations between PDMP without boundary and point processes are developed in the book of M. Jacobsen \cite{Jacobsen}.  Recently C. Cocozza-Thivent has deeply investigated relations between PDMP and Markov renewal theory \cite{Cocozza_2} and extended PDMP's  definition.  In all application areas most of interest quantities depend on the distribution of the process at each time, so it is essential of knowing how to compute these marginal distributions. The method proposed here consists in solving numerically equations which are fulfilled by the marginal distributions, namely generalized Kolmogorov equations. The characterization of the marginal distributions by these equations is studied in \cite{CEMR} for a PDMP without boundary. Finite volume schemes are proposed in
\cite{BEP, CE, CEM, EM, EMP, EMR, LMRZ}.

One studies the following class of PDMP with boundary.  The state space of the process is an open subset $F$ of $\xR^d$ and there exists a subset $\Gamma$ of the topological boundary of $F$ which will force the process to jump. The process $(X_t)_{t\ge0}$ is a jump stochastic process on $F$ whose trajectories are deterministic between the jump times. 

The deterministic trajectories are determined by a flow $\phi(x,t)$: if between $s$ and $t$ ($s < t$) the process does not reach the frontier and does not jump, then $X_{t} = \phi(X_s,t-s)$. Of course the flow has the ``Markov property'' $\phi(\phi(x,s),t)=\phi(x,s+t)$ as far as the boundary is not reached. Two kinds of jumps can occur. First there are stochastic jumps  from a position $x \in F$ with a jump rate $\lambda(x)$ and a jump distribution $Q(x, \xdif y)$. Second when the process reaches a point $x$ of the boundary $\Gamma$, it jumps inside $F$ with the distribution $q(x, \xdif y)$. Roughly speaking, these two kinds of jumps have different characters. The first ones occur at random times with probability density functions while the second ones occur at times with Dirac distributions. \medskip


We assume the following notations and hypotheses on the data, denoted by (H) in this paper. 
\begin{enumerate}
\item \label{hypE}  $d\in{\xN}^{\star}$ and  $\mathcal{P}(\xR^d)$ is the set of probability measures on
$\xR^d$ with borelian algebra, $\mathcal{P}(A)$ is the subset of $\mathcal{P}(\xR^d)$ with support in $A$ for any measurable subset $A\subset \xR^d$.
\item \label{hypflow} 
The flow  $\phi~:~\xR^{d}\times\xR_+\to \xR^{d}$ is assumed to be such that:

\begin{enumerate}
\item  \label{hypflowLip} $\phi~:~\xR^{d}\times\xR_+\to \xR^{d}$ is Lipschitz continuous with constant $L_\phi$,
\item \label{hypflowflow} $\phi(x,0) = x$ for all $x\in \xR^d$ and
\[
\forall x\in\xR^d,\forall t,s\in\xR_+, \qquad \phi(\phi(x,t),s) = \phi(x,t+s)
\]
\item \label{hypflowalpha} Let $F\subset\xR^d$ be a non empty open set and $G = {\xR}^d\setminus F$ the complementary of $F$ be  a non-empty closed set such that, for all $x\in \xR^{d}$, there exists $t\in {\xR}_{+}$ such that $ \phi(x,t)\in G$. We then define $\alpha~:~\xR^{d}\to {\xR}_{+}$ by
\[
 \alpha(x) = \inf\{t\ge 0, \phi(x,t)\in G\}.
\]
Note that, for all $x\in G$, $\alpha(x) = 0$, and that, for all
$x\in F$, since $\phi$ is continuous and $G$ is closed, $\alpha(x)>0$.  Note that the following property holds
\[
 \forall x\in F,\ \forall t\in (0,\alpha(x)), \alpha(\phi(x,t)) = \alpha(x) - t.
\]
 We assume that the function $\alpha$ is Lipschitz continuous with constant $L_\alpha$. 
\item We then denote $\Gamma = \{\phi(x,\alpha(x)),x\in F\}$. We have $\Gamma\subset G$. We cannot state
whether $\Gamma$ is open or closed.
\end{enumerate}

\item \label{hypLambda} The transition rate $\lambda$ is such that
$\lambda\in \mathrm{C}_{b}(F,{\xR}_{+})$, where $\mathrm{C}_{b}(F,{\xR}_{+})$ denotes the set of
continuous and bounded functions from $F$ to ${\xR}_{+}$. We denote by $\Lambda>0$ a bound  of $\lambda$. 

\item \label{hypQ} The transition probability $Q~:~F\to \mathcal{P}(F)$ (we then denote it by $x\mapsto Q(x,\xdif y)$) is such that:
\begin{enumerate}
\item \label{hypQfQ} there exists a function $f_Q~:~\xR_+ \to \xR_+$ such that
\[
f_Q(r) = \sup_{x\in F} \int_{\{y \in F : |y| \ge |x|+r\} }Q(x,\xdif y) \quad\text{ and }\quad\lim_{r\to\infty} f_Q(r) = 0.
\]

\item \label{hypQC0} for all $\xi\in \mathrm{C}_{b}(F,{\xR
})$, the function $x\rightarrow\int\xi(y) Q(x,\xdif y)$ is continuous from
$F$ to ${\xR}$.
\end{enumerate}

\item \label{hypq} The transition probability $q~:~ \Gamma \to \mathcal{P}(F)$  (we then denote it by $x\mapsto q(x,\xdif y)$) is such that:
 
\begin{enumerate}
\item \label{hypqfq} there exists a function $f_q~:~\xR_+ \to \xR_+$ such that
\[
f_q(r) = \sup_{x\in \Gamma} \int_{\{y \in F :  |y| \ge |x|+r\} }q(x,\xdif y)
 \quad\text{ and }\quad
\lim_{r\to\infty} f_q(r) = 0.
\]

\item \label{hypqC0} for all $\xi\in \mathrm{C}_{b}(F,{\xR
})$, the function $x\rightarrow\int\xi(y) q(x,\xdif y)$ is continuous from
$\Gamma$ to ${\xR}$.

\item\label{hypqexpB}  denoting $0 = \exp(-B\infty)$ for all $B>0$, 
 we assume that there exists $a_0\in (0,1)$ and $B_0>0$ such that
\[
\sup_{x\in\Gamma} \int_F e^{-B_0\alpha(y)}q(x,\xdif y) \le 1 - a_0.
\]

\end{enumerate}

\item \label{hyprho} We assume that $\rho_{\rm ini} \in \mathcal{P}(F)$ is given.
\end{enumerate}

These hypotheses are very natural because we only assume Lipschitz continuous regularity. However some controls (especially hypothesis (H.\ref{hypqexpB})) are technical and could certainly be modified to be more general. In general, the flow $\phi$ is more regular as solution of an ordinary differential equation of the form $\partial_{t}\phi(x,t) = v(\phi(x,t))$. If the flow $\phi$ has sufficient regularity, the function $\alpha$ is only the length of a trajectory and is well-defined. But the regularity of $\alpha$ crucially depends on the set $G$ (for instance there is difficult questions about the topology of this domain) and it is not clear that the Lipschitz continuity will be verified in many cases. However, there are many ways to choose the state space $F$ and the flow $\phi$ in order to obtain the same PDMP, thus there is certainly an other model (for instance with decomposition or duplication of the state space) with desired properties. The regularity in hypotheses (H.\ref{hypQC0}) and (H.\ref{hypqC0}) are very natural, and again the state space could certainly be modified in order to satisfy these hypotheses. In general $Q$ is absolutely continuous with respect to the Lebesgue measure such that for all $x\in F$ there exists $h_{x}$ a function on $F$ such that $\int \xi(y) Q(x, \xdif y) = \int \xi(y) h_{x}(y) \xdif y$ for all $\xi\in \mathrm{C}_{b}(F,{\xR})$, then the regularity is transferred to the density $h_{x}$. The hypothesis is now verified for a bounded continuous function $x \mapsto h_{x}$.

\bigskip

The distributions $\rho_t\in \mathcal{P}(F)$ of the process at time $t$ satisfy the generalized Kolmogorov equation
\begin{eqnarray} \label{eqn1}
\lefteqn{\int_{F} g(x,T) \rho_{T}(\xdif x) =
\displaystyle \int_F g(x,0) \rho_{\rm ini}(\xdif x) + \int_{F\times[0,T)} \partial_{t,\phi}g(x,t) \; \mu(\xdif x,\xdif t)}\\
& \displaystyle+ \int_{F\times[0,T)} \lambda(x) \left(\int_{F}  g(y,t) Q(x,\xdif y)  -  g(x,t)\right)  \mu(\xdif x,\xdif t) \nonumber\\
& \displaystyle+   \int_{\Gamma\times[0,T)}\left(\int_{F}  g(y,t)q(x,\xdif y) - g(x,t)\right) \sigma(\xdif x,\xdif t),\
&\forall g\in \mathcal{T}\nonumber
\end{eqnarray}
where $\mu(\xdif x,\xdif t) = \rho_t(\xdif x) \xdif t$. The measure $\sigma$, the test function space $\mathcal{T}$ and the operator $\partial_{t,\phi}$ are defined in Section \ref{S:uniqueness}. The measure $\sigma$ describes the average number of times that the trajectories reach some parts of the boundary. 


\section{Uniqueness}\label{S:uniqueness}

In the following, we propose a numerical scheme in order to compute numerically solutions of equation \eqref{eqn1}. A result of uniqueness is essential and needs an adapted test function space. This leads to the following definition.

Let us denote by $\mathrm{C}_b^c(\xR^d\times\xR_+)$ the set of all $g\in \mathrm{C}_b(\xR^d\times\xR_+)$ with compact support in time, meaning that there exists $T\in\xR_+$ (depending on $g$) such that $g(x,t)  = 0$ for all $t\ge T$.

For all $g \in \mathrm{C}_b^c(\xR^d\times\xR_+)$, let:
\[
\partial_{t,\phi} g(x,t) = \limsup_{n \to \infty} \frac{g(\phi(x,1/n), t+1/n) - g(x,t)}{1/n}
\]
when this limit is finite and $0$ if not.
The operator $\partial_{t,\phi}$ is called the derivation along the flow. If $\partial_{t}g(x,t) = \lim_{\epsilon \to 0} \frac{g(x,t+\epsilon) - g(x,t)}{\epsilon}$ and $\partial_{\phi}g(x,t) = \lim_{\epsilon \to 0} \frac{g(\phi(x,\epsilon),t) - g(x,t)}{\epsilon}$ exist, then for all $(x,t) \in \xR^{d}\times\xR_+$ and $0 \le t < \alpha(x)$ the function $t \mapsto g(\phi(x,t),t)$ is differentiable and $\partial_{t,\phi}g(x,t) = \partial_{t}g(x,t) + \partial_{\phi}g(x,t)$.

\begin{dfntn}
We denote by $\mathcal{T}$ the set of all functions $g\in \mathrm{C}_b^c(\xR^d\times\xR_+)$  such that there exists
$I,J\in \mathrm{C}_b^c(\xR^d\times\xR_+)$ with
\[
\forall (x,t)\in \xR^{d}\times\xR_+,\  g(x,t) = J(\phi(x,\alpha(x)),t+\alpha(x)) - \int_0^{\alpha(x)} I(\phi(x,s),t+s)\xdif s.
\]
We then denote $g = \mathbb{T}(I,J)$.\end{dfntn}
We will use $\mathcal{T}$ as test function space in the proof of uniqueness. If $g = \mathbb{T}(I,J) \in \mathcal{T}$, it is easy to verify that  for $x\in G$ we have $g(x,t) = J(x,t)$ and that for $x\in F$ and $0 < \epsilon < \alpha(x)$ we have
\[
 g(\phi(x,\epsilon),t+\epsilon)-g(x,t) =  \int_0^{\epsilon} I(\phi(x,s),t+s)\xdif s.
\]
Then, for all $g = \mathbb{T}(I,J)\in \mathcal{T}$ and for all $x\in F$ and $t\in \xR_+$, we have
$I = \partial_{t,\phi}g$. Hence, for given $I,J,\widetilde I,\widetilde J \in \mathrm{C}_b^c(\xR^d\times\xR_+)$, such that $g = \mathbb{T}(I,J) = \mathbb{T}(\widetilde I,\widetilde J)$, we get $I = \widetilde I$ on $F$. We also have $J = \widetilde J = g$ on $G$.

Now let us define precisely the problem for which uniqueness is proven.
\begin{dfntn} \label{probleme}
We say that non negative measures $\mu$ and $\sigma$ such that
$\forall \ T\in \xR_+,\ \mu(F \times [0,T])<\infty, \ \sigma(\Gamma\times [0,T])<\infty,\ \mu((\xR^d\setminus F) \times \xR_+)=0\hbox{ and }\sigma((\xR^d\setminus\Gamma) \times \xR_+)=0$
are solutions of Problem {\bf{P}}  if for all $g = \mathbb{T}(I,J) \in \mathcal{T}$
\begin{eqnarray}\label{eqn2bis}
&\displaystyle 0 = \int_F g(x,0) \rho_{\rm ini}(\xdif x) + \int_{F\times\xR_+} I(x,t)\mu(\xdif x,\xdif t)\nonumber\\
&\displaystyle + \int_{F\times\xR_+ }\lambda(x) \left(\int_{F} g(y,t) Q(x,\xdif y)  - g(x,t)\right) \mu(\xdif x, \xdif t)\label{eq:musigma}\\
&\displaystyle +  \int_{\Gamma\times\xR_+}\left(\int_{F}  g(y,t)q(x,\xdif y) - g(x,t)\right) \sigma(\xdif x,\xdif t). \nonumber
\end{eqnarray}
\end{dfntn}

We have the following uniqueness result.

\begin{thrm}\label{uniqueness}
Under hypotheses (H), there exists at most a unique couple $(\mu, \sigma)$ solution of Problem {\bf{P}}.
\end{thrm}

\begin{proof}
Suppose there exist two solutions $(\tilde{\mu},\tilde{\sigma})$ and $(\hat{\mu},\hat{\sigma})$ to  Problem {\bf{P}}. Denote $(\bar{\mu},\bar{\sigma})$ the measures such that $\bar{\mu} = \tilde{\mu} - \hat{\mu}$ and $\bar{\sigma} = \tilde{\sigma} - \hat{\sigma}$.
Then for all $g = \mathbb{T}(I,J) \in\mathcal{T}$, we have
\begin{eqnarray}
\lefteqn{0 = \int_{F\times\xR_+ } I(x,t)\bar{\mu}(\xdif x,\xdif t)}\nonumber\\
&&+ \int_{F\times\xR_+}\lambda(x) \left(\int_{F} g(y,t) Q(x,\xdif y)  - g(x,t)\right) \bar{\mu}(\xdif x, \xdif t) \nonumber\\
&&+ \int_{\Gamma\times\xR_+}\left(\int_{F}  g(x,t)q(z,\xdif x)- g(z,t) \right)\bar{\sigma}(\xdif z,\xdif t). \nonumber
\end{eqnarray}

Let $\overline I, \overline J \in \mathrm{C}_b^c(\xR^d\times\xR_+)$. Using Lemma \ref{op_inversion}, we can find $I,J\in \mathrm{C}_b^c(\xR^d\times\xR_+)$
such that $g =\mathbb{T}(I,J)$ verifies
\begin{eqnarray*}
&\displaystyle \forall (x,t) \in \xR^{d}\times \xR_{+},\ \overline I(x,t) = I(x,t) + \lambda(x) \left(\int_{F} g(y,t) Q(x,\xdif y)  - g(x,t)\right),\\
\text{and}& \displaystyle \forall (z,t) \in \xR^{d}\times \xR_{+},\ \overline J(z,t) = \int_{F}  g(x,t)q(z,\xdif x)- J(z,t).
\end{eqnarray*}

Thus the measures $\bar{\mu}$ and $\bar{\sigma}$ verify
\[
\int_{F\times\xR_+}  \overline I(x,t) \bar{\mu}(\xdif x, \xdif t) + \int_{\Gamma\times\xR_+} \overline J(z,t)\bar{\sigma}(\xdif z,\xdif t) = 0.
\]
Since this equality is verified for all $\bar{I},\bar{J} \in \mathrm{C}_b^c(\xR^d\times\xR_+)$, it proves that the measures $\bar{\mu}$ and $\bar{\sigma}$ vanish.\end{proof}
The previous proof uses the following lemma which is an original result. In fact, this is the core of the problem and it has given to the author a good representation for the space $\mathcal{T}$ of test functions. Moreover many hypotheses are used in this lemma and the reader can understand the choice of hypothesis by reading this proof, in particular it explains the form of hypothesis (H.\ref{hypqexpB}).

\begin{lmm}[Operator's inversion]\label{op_inversion} Under hypotheses (H), let $\overline I, \overline J \in \mathrm{C}_b^c(\xR^d\times\xR_+)$. Then
there exists  $I,J  \in \mathrm{C}_b^c(\xR^d\times\xR_+)$ such that, setting $g =\mathbb{T}(I,J)$, we have
\begin{equation}\label{EqI}
\forall (x,t) \in \xR^{d}\times \xR_{+},\ \overline I(x,t) = I(x,t) + \lambda(x) \left(\int_{F} g(y,t) Q(x,\xdif y)  - g(x,t)\right),
\end{equation}
and 
\begin{equation}\label{EqJ}
\forall (z,t) \in \xR^{d}\times \xR_{+},\ \overline J(z,t) = \int_{F}  g(x,t)q(z,\xdif x)- J(z,t).
\end{equation}
\end{lmm}
\begin{proof} There exists $T>0$ such that $\overline{I}(x,t)=\overline{J}(x,t)=0$ for all $x\in\xR^{d}$ and $t\geq T$. Let
\[
\mathrm{C}_b^T = \{ g\in \mathrm{C}_b(\xR^d\times\xR_+),\ \forall (x,t)\in \xR^d\times [T,+\infty[, g(x,t)  = 0\}.
\]
Remark that it is clear that for all $I,J\in\mathrm{C}_b^T$, $\mathbb{T}(I,J) \in \mathrm{C}_b^T$. We define $\alpha^{T}(x) = \min(\alpha(x),T)$ for all $x\in \xR^d$ and a norm on $\mathrm{C}_b^T$ by
\begin{equation*}
\forall g \in \mathrm{C}_b^T, \quad \| g \|_{A,B} :=\sup_{(x,t) \in \xR^d\times [0,T]} \exp(B\alpha^T(x)+ At)\ \ | g(x,t)|,
\end{equation*}
for given $A,B>0$ chosen later. $\mathrm{C}_b^T$ is a Banach space with this norm.

Then we define a sequence of functions $(I^{n},J^{n} )_{n\in \xN}$ such that $I^{0} = J^{0} = 0$ and 
for all $n\in\xN$
\begin{eqnarray*}
& \displaystyle \forall (x,t) \in \xR^{d}\times \xR_{+},\ I^{n+1}(x,t) = \overline  I(x,t) - \lambda(x) \left(\int_{F} \mathbb{T}(I^{n},J^{n})(y,t) Q(x,\xdif y)  - \mathbb{T}(I^{n},J^{n})(x,t)\right),\\
& \displaystyle \forall (z,t) \in \xR^{d}\times \xR_{+},\  J^{n+1}(z,t) =  \int_{F}  \mathbb{T}(I^{n},J^{n})(x,t)q(z,\xdif x) - \overline  J(z,t).
\end{eqnarray*}
This define a sequence $(\mathbb{T}(I^n,J^n))_{n\in\xN}$ in $\mathrm{C}_b^T$ where $\mathbb{T}(I^{n+1},J^{n+1}) = \Psi(\mathbb{T}(I^n,J^n))$ for an obvious operator $\Psi$ . Let us prove that this sequence converges. For this, it is sufficient to prove that the operator $\Psi$ is really well defined and there exists $k\in (0,1)$ such that for all $(I,J) \in (\mathrm{C}_b^T)^2$ and $(I',J') \in (\mathrm{C}_b^T)^2$ we have
\[
 \| \Psi(\mathbb{T}(I,J)) - \Psi(\mathbb{T}(I',J')) \|_{A,B} \le k  \| \mathbb{T}(I,J)-\mathbb{T}(I',J') \|_{A,B}.
\]
Then, using the definition of $(I^{n+1},J^{n+1} )$, we get 
that the sequence $(I^{n},J^{n} )_{n\in \xN}$ is convergent as well. The limit satisfy \eqref{EqI} and \eqref{EqJ}.\medskip

Let us set  $f = \mathbb{T}( I-I', J-J')$, let $(x,t)\in \xR^d\times[0,T]$. We have
\[
\left(\Psi(\mathbb{T}(I,J)) - \Psi(\mathbb{T}(I',J'))\right)(x,t) = T_1 - T_2 + T_3,
\]
with
\begin{equation*}
\displaystyle T_1 = \int_0^{\alpha^T(x)} \lambda(\phi(x,s)) \int_{F} f(y,t+s) Q(\phi(x,s),\xdif y) \xdif s,
\end{equation*}
\begin{equation*}
\displaystyle T_2 = \int_0^{\alpha^T(x)} \lambda(\phi(x,s))  f(\phi(x,s),t+s)\xdif s,\quad\text{and}\quad  T_3 = \int_{F}  f(y,t+\alpha(x)) q(\phi(x,\alpha(x)),\xdif y).
\end{equation*}
We then have
\begin{eqnarray*}
|T_1|&\leq & \Big|\int_{0}^{\alpha^T(x)}\exp(-A(t+s)) \lambda(\phi(x,s))\\
&&\qquad\times\int_{F} \exp(A(t+s))\exp(-B\alpha^{T}(y))\exp(B\alpha^{T}(y))f(y,t+s)Q(\phi(x,s),\xdif y) \xdif s\Big|\\
&\leq &\int_{0}^{\alpha^T(x)}\exp(-A(t+s)) |\lambda(\phi(x,s))|\int_{F} \exp(-B\alpha^{T}(y))\|f\|_{A,B}Q(\phi(x,s),\xdif y) \xdif s\\
&\leq &\|f\|_{A,B} \int_{0}^{\alpha^T(x)}\exp(-A(t+s)) |\lambda(\phi(x,s))| \xdif s.
\end{eqnarray*}
Therefore we get
\begin{eqnarray*}
\exp(At+B\alpha^T(x))\left|T_1\right|&\leq& \|f\|_{A,B}\int_{0}^{\alpha^T(x)}\exp(-As)|\lambda(\phi(x,s))|\exp(B\alpha^T(x)) \xdif s\\
&\leq&  \frac{\|f\|_{A,B} \Lambda \exp(BT)}{A}.
\end{eqnarray*}
In the same way, assuming $A> B$ and using $\alpha(\phi(x,s)) = \alpha(x) - s$ which leads to $\alpha^{T}(\phi(x,s)) \ge \alpha^T(x) - s$,
\begin{eqnarray*}
\exp(At+B\alpha^T(x)) \left |T_2 \right |
&\leq& \Lambda\|f\|_{A,B} \frac{1}{A-B}.
\end{eqnarray*}
Finally, for $A>B$, using $f(y,t+\alpha(x)) = f(y,t+\alpha^{T}(x))$ for all $(y,t)\in\xR^{d}\times\xR_{+}$ for all $x\in\xR^{d}$, we have 
\begin{eqnarray*}
\exp(At+B\alpha^T(x))\left|T_3 \right |
&\leq&\|f\|_{A,B} \left(\int_{F}\exp(-B \alpha(y)) q(\phi(x,\alpha(x)),\xdif y) + \exp(-B T) \right).
\end{eqnarray*}

Finally we need to choose sufficiently large constants $A$ and $B$ with $A>B$ such that 
\[
\sup_{x\in \xR^d}\left(\frac{\Lambda\exp(BT)}{A} +\frac{\Lambda}{A-B} +\int_{F}\exp(-B \alpha(y))q(\phi(x,\alpha(x)),\xdif y) + \exp(-BT) \right) \le k < 1.
\]
Thanks to Hypothesis (H.\ref{hypqexpB}), setting first large $B$ and next large $A$, we obtain the result for $k=1-\frac{a_{0}}{2}<1$.

So it only remains to prove that the operator $\Psi$ is well defined such that it effectively maps $\mathcal{T}\cap \mathrm{C}_{b}^{T}$ to $\mathcal{T}\cap \mathrm{C}_{b}^{T}$. Remark that it is sufficient to prove that the following terms $(U_{i})_{i=1..5}$ are in $\mathrm{C}_{b}^{T}$.
\begin{eqnarray*}
&U_{1}:=\displaystyle \int_{0}^{\alpha(x)} \overline{I}(\phi(x,s),t+s)\xdif s, \quad U_{2}:= \displaystyle \int_{0}^{\alpha(x)}
\lambda(\phi(x,s)) \int_{F} \mathbb{T}(I,J)(y,t+s) Q(\phi(x,s),\xdif y)\xdif s,\\
&U_{3}:=\displaystyle \int_{0}^{\alpha(x)}
\lambda(\phi(x,s))\ \mathbb{T}(I,J)(\phi(x,s),t+s)\xdif s,\\
&U_{4}:= \displaystyle \int_{F} \mathbb{T}(I,J)(y,\alpha(x)+t) q(\phi(x,\alpha(x)),\xdif y)\quad \text{ and }\quad U_{5}:=\overline{J}(\phi(x,\alpha(x)),\alpha(x)+t).
\end{eqnarray*}
First remark that all terms are bounded since
\[
|U_{1}| \leq T\|\overline{I}\|_{\infty}, \quad |U_{2}|+|U_{3}| \leq 2T \Lambda \|\mathbb{T}(I,J)\|_{\infty}, \quad |U_{4}| \leq  \|\mathbb{T}(I,J)\|_{\infty} \quad \text{ and } \quad |U_{5}| \leq \|\overline{J}\|_{\infty}.
\]
Moreover, for all $t > T$, all terms vanish. So it only remains to prove the continuity. Since $\overline{I}$, $\overline{J}$, $\phi$ and $\alpha$ are continuous, $U_{1}$ and $U_{5}$ are clearly continuous. $U_{3}$
is clearly continuous since $\lambda$ and $\mathbb{T}(I,J)$ are bounded continuous functions. Finally $U_{2}$ and $U_{4}$ are continuous too: it is relatively straightforward using Hypotheses (H.\ref{hypQfQ}), (H.\ref{hypQC0}), (H.\ref{hypqfq}) and (H.\ref{hypqC0}). \end{proof}


\section{A finite volume scheme}\label{FVscheme}

We now come to the presentation of the
finite volume scheme, which has been used in \cite{EM} for the approximation 
of a benchmark problem. Such schemes are not classically used in the framework
of probabilistic studies, since they have mainly be developed by engineers,
in order to approximate the solutions of balance equations. We then consider that \eqref{eqn1} can
be viewed as balance equations describing the conservation of probability. 
We then introduce the finite volume discretization by the following notations and definitions.

\begin{enumerate}
\item We define a reference measure, denoted by $\xdif x$ or $\xdif y$, on $F$, with respect to all borelian sets of $\xR^d$ restricted to $F$.

\item An admissible mesh $\mathcal{M}$ of $F$ is a countable partition of $F$, therefore such that
$\cup_{K\in\mathcal{M}}K=F$ and $\forall(K,L)\in\mathcal{M}^{2},K\neq L\Rightarrow K\cap L=\emptyset$.

\item $\forall K\in\mathcal{M},\ |K|:=\int_{K}\xdif x>0$.

\item $\sup\limits_{K\in\mathcal{M}}\mathrm{diam}(K)<+\infty$ where
$\mathrm{diam}(K)=\sup\limits_{\left\{  \left(  x,y\right)  \in K^{2}\right\}
}\left|  x-y\right|  $. We then set $h :=\sup\limits_{K\in\mathcal{M}}\mathrm{diam}(K)$.

\item $\tau>0$ and $\dt>0$  are given values, and we denote by $\mathcal{D} = (\mathcal{M},\delta\!t,\tau)$.
\end{enumerate}

The value $\tau>0$, aimed to tend to 0, is used for the definition, for all $K\in\mathcal{M}$ and $L\in\mathcal{M}$, of the flux of
probability mass from $K$ to $L$ by
\begin{equation}
v_{KL}=\frac{1}{\tau}|\{x\in K: \alpha(x)>\tau\hbox{ and } \phi(x,\tau)\in L\}|.
\label{defvkl}
\end{equation}
And $\forall(K,L)\in\mathcal{M}\times\mathcal{M}$ we define
\begin{equation}
\displaystyle\lambda_{KL} =\int_{K} \lambda(x)\int_L Q(x,\mathrm{d}y) \mathrm{d}x, \qquad
\displaystyle \lambda_{K}  =\int_{K}\lambda(x)\mathrm{d}x=\sum_{L\in\mathcal{M}}\lambda_{KL},
\label{deflambdakl}
\end{equation}
\begin{equation}
\displaystyle q_{K} =\frac{1}{\tau}|\{x\in K: \alpha(x)\le {\tau}\}|,\qquad\text{and}\qquad  \displaystyle q_{KL} =\frac{1}{{\tau}} \int_{\{x\in K: \alpha(x)\le {\tau}\}} \int_L q(\phi(x,\alpha(x)),\mathrm{d}y) \mathrm{d}x.
\label{defqkl}
\end{equation}
We may now define a family $\left(  p_{n}^{(K)}\right)_{n\in\xN,K\in\mathcal{M}}$
 of real values thanks to the following finite volume scheme, which is time implicit.
\begin{align}
&  |K|\frac{p_{n+1}^{(K)}-p_{n}^{(K)}}{\delta\!t}+\sum_{L\in\mathcal{M}}\left(  v_{KL}p_{n+1}^{(K)} -v_{LK}p_{n+1}^{(L)}\right)  \nonumber\\
&  + (\lambda_{K}+q_{K}) p_{n+1}^{(K)}- \sum_{L\in\mathcal{M}}p_{n+1}^{(L)} (\lambda_{LK}  + q_{LK}) = 0,& 
\forall K  \in\mathcal{M},\ \forall n\in{\xN},
\label{schnum}
\end{align}
with the initial condition
\begin{equation}
|K|\,p_{0}^{(K)}=\int_{K}\rho_{\rm ini}(\mathrm{d}x),\qquad \forall K\in\mathcal{M}.
\label{schnumini}
\end{equation}
Let us remark that the following property holds:
\begin{equation}
{\tau}\left(\sum_{L\in\mathcal{M}} v_{KL} + q_{K}\right)= |K|,\qquad \forall K\in\mathcal{M}.
\label{propvkl}
\end{equation}
Therefore, scheme \eqref{schnum} may be rewritten
\begin{equation}
\left(\left(1 + \frac {\delta\!t}{\tau}\right)|K|+{\delta\!t}\lambda_{K}\right) p_{n+1}^{(K)} -{\delta\!t}\sum_{L\in\mathcal{M}}p_{n+1}^{(L)} (v_{LK} + \lambda_{LK}  + q_{LK})=|K|\,p_{n}^{(K)},\quad
\forall K \in\mathcal{M},\ \forall n\in{\xN}.
\label{schnumis}
\end{equation}
We then define the approximation $\mu_{\mathcal{D}}(\xdif x,\xdif t)$ (resp. $\sigma_{\mathcal{D}}(\xdif x,\xdif t)$) of the measure $\mu(\xdif x,\xdif t)$ on $F\times\xR_+$ (resp. $\sigma(\xdif x,\xdif t)$ on $\Gamma\times\xR_+$) by
\begin{equation*}
\int_{\xR^d\times\xR_+} f(x,t) \mu_{\mathcal{D}}(\xdif x,\xdif t)=  \sum_{n\in \xN}{\delta\!t}  \sum_{K\in\mathcal{M}} p_{n+1}^{(K)} \int_K f(x,n\delta\!t)\mathrm{d}x,
\end{equation*}
for all bounded continuous function $f\in\mathrm{C}_b^c$, and
\begin{equation*}
\int_{\xR^d\times\xR_+} f(x,t) \sigma_{\mathcal{D}}(\xdif x,\xdif t) =  \sum_{n\in \xN}{\delta\!t} \sum_{K\in\mathcal{M}}  p_{n+1}^{(K)}\frac{1}{{\tau}} \int_{\{x\in K: \alpha(x)\le {\tau}\}} f(\phi(x,\alpha(x)),n\delta\!t)\mathrm{d}x,
\end{equation*}
for all bounded continuous function $f\in\mathrm{C}_b^c$ in such a way that $\mu_{\mathcal{D}}$ and $\sigma_{\mathcal{D}}$ are expected to converge to $\mu$ and $\sigma$.

The specifications of this scheme, depending on the mesh and two parameters, $\tau$ and $\dt$, are resulting from the following observations:
\begin{enumerate}
\item An explicit version of the scheme could be defined, following the ideas of \cite{CEM}.
But the main drawback of an explicit scheme is that it cannot provide, in the general case, an approximation of the asymptotic marginal distributions at large times, with an acceptable computing cost.
\item An implicit scheme has been provided in \cite{EMP}. But in this scheme the considered flow is much more regular (it is assumed to be the solution of an EDO).
\item An explicit scheme is introduced in \cite{BEP} for Lipschitz flows. The present scheme uses the value $\tau>0$ in the same way as the time step is used in  \cite{BEP} (where the convergence of the scheme is proved for general Lipschitz flow  $\phi$ in the case $\dt\to 0$ and $h/ \dt \to 0$).
\item In \cite{EM}, the asymptotic states at large times have been obtained letting $\dt\to\infty$ in an implicit scheme. Hence it is interesting to use the parameters $h$ and $\tau$ which can tend to $0$ independently of $\dt$.
\end{enumerate}

\begin{rmrk}
Since \eqref{schnum} is an infinite linear system, the existence and uniqueness of a positive solution must be first addressed. From a numerical point of view, this infinite linear system should be treated with caution, but we often have a system with finite number of unknowns.
\end{rmrk}

\begin{lmm} [Existence of solution]\label{lemexistsol}
Under Hypotheses (H) and numerical hypotheses, there exists one and only one solution $(p_{n}^{(K)})_{K \in\mathcal{M},n\in{\xN}}$ to Scheme \eqref{defvkl}, \eqref{defqkl}, \eqref{deflambdakl}, \eqref{schnum}, \eqref{schnumini} which satisfies:
\begin{eqnarray}
&p_{n}^{(K)} \ge 0,& \forall K \in\mathcal{M},\ \forall n\in{\xN},\label{positive}\\
&\sum_{K \in\mathcal{M}} |K| \, p_{n}^{(K)} = 1,& \forall n\in{\xN}\label{massproba}.
\end{eqnarray}
\end{lmm}
\begin{proof}
Let us first show the existence of a solution to the scheme. We consider
the values $p_{(k)}^{(K)}$ defined, for given $n\in{\xN}$ and $(p_{n}^{(K)})_{K \in\mathcal{M}}$ such that \eqref{positive}-\eqref{massproba}, by:
\begin{equation}
\left\{\begin{array}{lr}
p_{(0)}^{(K)} = p_{n}^{(K)},&\forall K \in\mathcal{M},\\
\\
\displaystyle\left(\left(1 + \frac {\delta\!t}{\tau}\right)|K|+{\delta\!t}\lambda_{K}\right) p_{(k+1)}^{(K)} = {\delta\!t}\sum_{L\in\mathcal{M}}p_{(k)}^{(L)} (v_{LK} + \lambda_{LK}  + q_{LK})+ |K|\,p_{n}^{(K)},&\forall K \in\mathcal{M},\ \forall k\in{\xN}.
\label{schnumfixpoint}
\end{array}\right.
\end{equation}
Denoting, for $k\in{\xN}$ and $K  \in\mathcal{M}$, $\widehat p_{(k+1)}^{(K)}= p_{(k+1)}^{(K)} - p_{(k)}^{(K)}$, we have
\begin{equation*}
\left\{\begin{array}{lr}
\displaystyle\left(\left(1 + \frac {\delta\!t}{\tau}\right)|K|+{\delta\!t}\lambda_{K}\right)\widehat p_{(1)}^{(K)} = {\delta\!t}\sum_{L\in\mathcal{M}} p_{n}^{(L)} (v_{LK} + \lambda_{LK}  + q_{LK}) - \left(\frac {\delta\!t}{\tau}|K|+{\delta\!t}\lambda_{K}\right) p_{n}^{(K)},&\forall K \in\mathcal{M},\\
\displaystyle\left(\left(1 + \frac {\delta\!t}{\tau}\right)|K|+{\delta\!t}\lambda_{K}\right)\widehat p_{(k+1)}^{(K)} = {\delta\!t}\sum_{L\in\mathcal{M}}\widehat p_{(k)}^{(L)} (v_{LK} + \lambda_{LK}  + q_{LK}),& \!\!\!\!\!\forall K  \in\mathcal{M},\forall k\in{\xN}^\star.
\end{array}\right.
\end{equation*}
We notice that, thanks to  \eqref{propvkl} and \eqref{massproba}, we have
\[
 \sum_{K \in\mathcal{M}} \left(\left(1 + \frac {\delta\!t}{\tau}\right)|K|+{\delta\!t}\lambda_{K}\right)\left|\widehat p_{(1)}^{(K)}\right| \le 2 \left(\frac {\delta\!t}{\tau} + \Lambda{\delta\!t}\right).
\]
and
\begin{align*}
&\sum_{K \in\mathcal{M}}\left(\left(1 + \frac {\delta\!t}{\tau}\right)|K|+{\delta\!t}\lambda_{K}\right) \left|\widehat p_{(k+1)}^{(K)}\right| \le \sum_{L\in\mathcal{M}}\left( \frac {\delta\!t}{\tau}|L|+{\delta\!t}\lambda_{L}\right) \left|\widehat p_{(k)}^{(L)}\right|,& \forall k\in\xN^\star.
\end{align*}
Then, by induction, we get that the value $u_k = \sum_{L\in\mathcal{M}}\left( \frac {\delta\!t}{\tau}|L|+{\delta\!t}\lambda_{L}\right) |\widehat p_{(k)}^{(L)}|$ is positive and nonincreasing
with respect to $k\in \xN^\star$. We then deduce that the sequence $(u_k)_{k\in \xN^\star}$ is convergent. Writing
\[
 \sum_{K \in\mathcal{M}} |K| \left|\widehat p_{(k+1)}^{(K)}\right| \le u_{k} - u_{k+1},
\]
we deduce that $\widehat p_{(k+1)}^{(K)}$ is the general term of an absolutely convergent series. Therefore the sequence $\left(p_{(k+1)}^{(K)}\right)_{k\in  \xN}$ satisfying
\[
 p_{(k+1)}^{(K)} = p_{n}^{(K)} + \sum_{m=0}^k \widehat p_{(m+1)}^{(K)},
\]
converges to a value denoted by $p_{n+1}^{(K)}$. Passing to the limit in \eqref{schnumfixpoint}, we obtain that these values are solution to the scheme. Moreover,
they satisfy that
\[
 \sum_{K \in\mathcal{M}}|K| p_{(k)}^{(K)} \le 1 + u_1 - u_{k} \le 1 + u_1.
\]
Therefore we can sum \eqref{schnum} on $K\in \mathcal{M}$, obtaining by induction that \eqref{massproba} holds for $n+1$.

\begin{rmrk}
In the case of an infinite state space, we can consider a bounded subdomain $\tilde{F}$ to do the computation until the density reaches the boundary of $\tilde{F}$. Next enlarge the subdomain $\tilde{F}$ and continue the computation. However, since the scheme is implicit, it could arrived after very few iteration.
\end{rmrk}
\begin{rmrk}
The matrix representing the linear system is potentially completely full because many coefficients $\lambda_{KL}$ are strictly positive and there is no simple way to reduce the complexity of this system inversion. But we only need to compute a unique LU decomposition to do multiple iterations of the scheme, because it is time independent. Moreover the computation of the coefficients of the matrix $v_{KL}$, $\lambda_{KL}$ and $q_{KL}$ is not an easy problem, because there is numerical difficulties depending on the geometry of the mesh or the quadrature rule used for the approximation of integrals.
\end{rmrk}

\medskip

Let us now turn to the uniqueness of this solution. Denoting by $\widehat p^{(K)}$ the difference between two solutions of 
 \eqref{schnumis}, and using that the two solutions satisfy  \eqref{positive}-\eqref{massproba},  one may write 
\begin{align*}
&  \left(\left(1 + \frac {\delta\!t}{\tau}\right)|K|+{\delta\!t}\lambda_{K}\right) |\widehat p^{(K)}| \le {\delta\!t}\sum_{L\in\mathcal{M}}|\widehat p^{(L)}| (v_{LK} + \lambda_{LK}  + q_{LK}),\ \forall K  \in\mathcal{M},
\end{align*}
which provides, summing on  $K  \in\mathcal{M}$, 
\begin{align*}
&  \sum_{K \in\mathcal{M}}\left(\left(1 + \frac {\delta\!t}{\tau}\right)|K|+{\delta\!t}\lambda_{K}\right) |\widehat p^{(K)}|  \le\sum_{L \in\mathcal{M}}\left(\frac {\delta\!t}{\tau}|L|+{\delta\!t}\lambda_{L}\right) |\widehat p^{(L)}|,
\end{align*}
and therefore $\sum_{K \in\mathcal{M}}|K| |\widehat p^{(K)}| = 0$. Hence the uniqueness proof.\end{proof}

\begin{rmrk} Since we have proved that the scheme is well-posed, some numerical tests are described in the Section \ref{conclusion}. But the convergence of the scheme needs some additional properties (whose proofs are not detailed here):
\begin{itemize}
\item finiteness of the mass $\sigma_{\mathcal{D}}(\Gamma\times[0,T])$,
\item tightness of the family of probability measures $(\mu_{\mathcal{D}})_{\mathcal{D}\in\mathcal{F}}$ given a family $\mathcal{F}$ of discretizations $\mathcal{D}$ verifying some restrictions like $h/\tau\leq1$,
\item tightness of the family of probability measures $(\sigma_{\mathcal{D}})_{\mathcal{D}\in\mathcal{F}}$ given the previous family $\mathcal{F}$.
\end{itemize}
\end{rmrk}
 
\section{Numerical test and conclusion}\label{conclusion}

Now we present different models in dimension $d=1$ to illustrate the efficiency of the numerical method for describing the evolution of the density of a PDMP with (or without) boundary. We are interested in the Transmission Control Protocol (TCP) window size process (described in \cite{DGR}). The first model is named the infinite maximum window size (TCP-I), the second model is a finite maximum window size (TCP-F) and the third model is a finite maximum window size with jump (TCP-FJ). In our framework, let the state space\footnote{Actually the state space should be a subset of $\xR_{+}$. But in our framework $F$ is open so we choose $(-\infty,X)$ in order to have $0\in F$. Of course the definition of $\phi$, $\alpha$, $\lambda$, $Q$ are only useful on $\xR_{+}$. Other tricks can be used to completely fulfills all the Hypotheses (H) and we do not give details here.} $F$ be the interval $(-\infty,X)$ for given $X>0$. Its complementary set $G$ is $[X,+\infty)$. The flow is given by $\phi(x,t) = \min(x+t,X)$ for all $(x,t)\in \xR\times\xR_{+}$. It is clear that $\alpha(x) = X-x$ for all $x\in F$ and $\alpha(x)=0$ for all $x\in G$, so $\Gamma = \{X\}$. The jump rate $\lambda$ is an increasing function of the position, precisely $\lambda : x \mapsto x^{+}$. Thus,``the higher a component is, the sooner it will jump''. $Q(x,\xdif y)$ is a Dirac measure on point $x/2$ for all $x\in F$. These assumptions are common to the three models. For the TCP-I model, we take $X$ as large as possible such that it does not perturb the dynamics. For the TCP-F model, we take a finite $X$ and give no particular behavior at $\Gamma$. Finally for the TCP-FJ model, we take a finite $X$ and force the process to jump. Precisely, at point $\{X\}$ the process jumps with probability $p $ at point $\{X/2\}$ and with probability $1-p$ uniformly on $(0,X)$ so 
\[
q(0,\xdif y) = p \ \delta_{X/2}(y) + \frac{1-p}{X}\  \chi_{(0,X)}(y) \xdif y,
\]
with obvious notations. Finally we choose $\rho_{\rm ini}$ as the Dirac distribution at point $\{0\}$
 (the computer starts).  The distributions $\rho_T$ of this process at time $T$ satisfy the generalized Kolmogorov equation
\begin{eqnarray*}
\lefteqn{\int_{(-\infty,X)} g(x,T) \rho_{T}(\xdif x) =
\displaystyle g(0,0) + \int_{(-\infty,X)\times[0,T)} (\partial_{t}+\partial_{x})g(x,t) \; \mu(\xdif x,\xdif t)}\\
& \displaystyle+ \int_{(-\infty,X)\times[0,T)}x \left(g(x/2,t) -  g(x,t)\right)  \mu(\xdif x,\xdif t)\\
& \displaystyle+ \int_{[0,T)}\left(p\ g(X/2,t) +(1-p) \int_{(0,X)}  g(y,t)\xdif y - g(X,t)\right) \sigma(\{X\},\xdif t),\
&\forall g\in \mathcal{T}.
\end{eqnarray*}

\begin{figure}[!htp]
\IGS{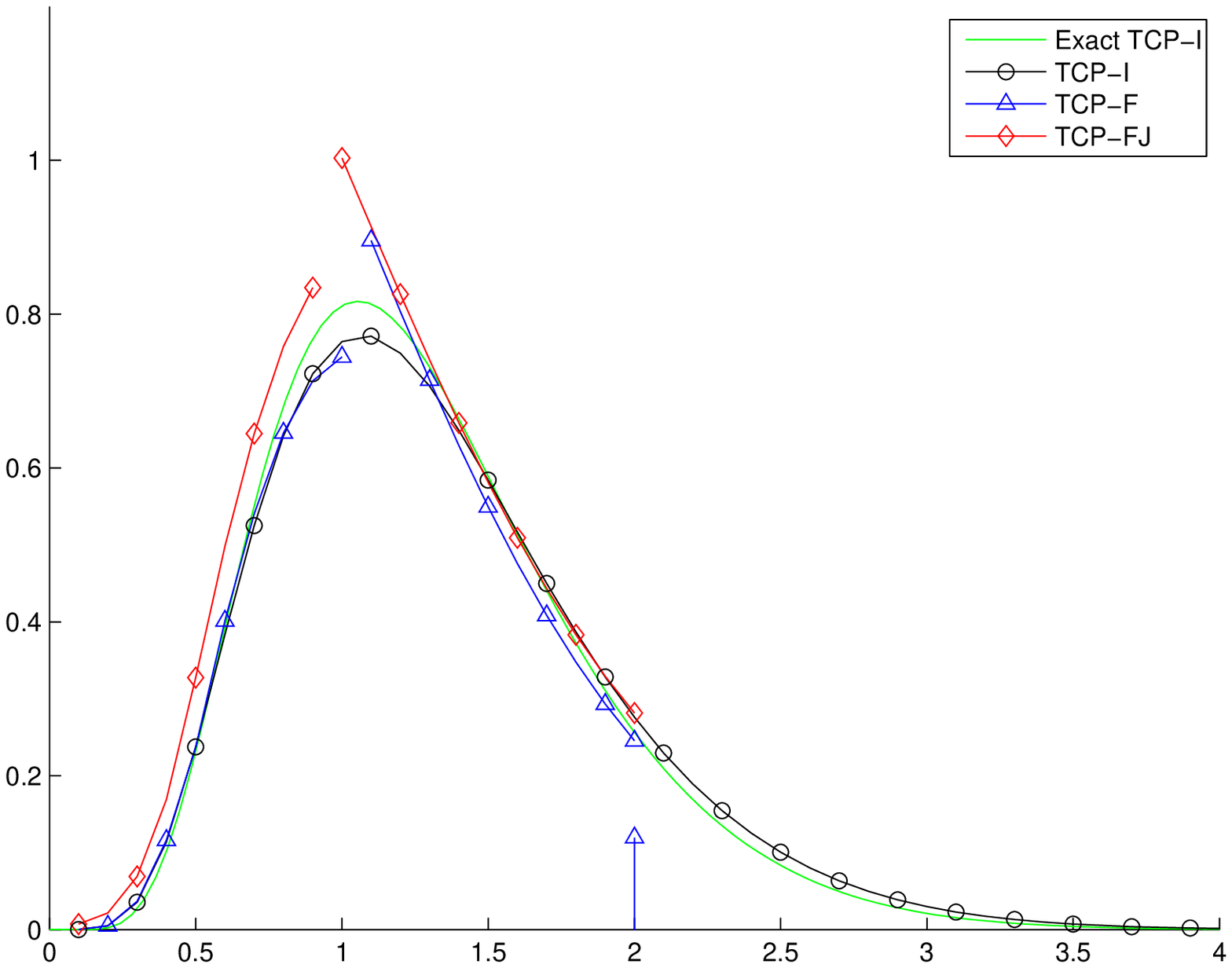}
\quad
\IGS{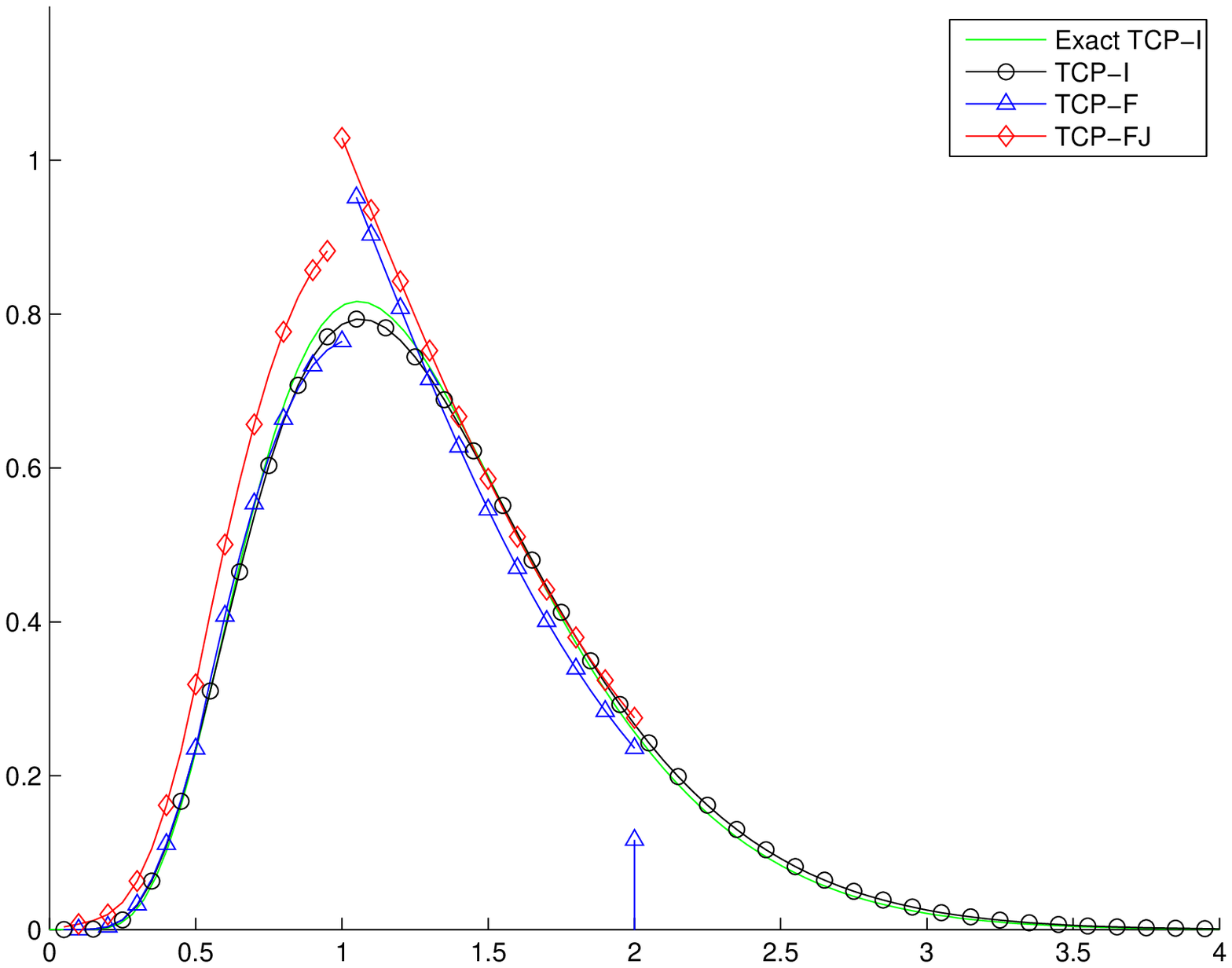}
\\
\IGS{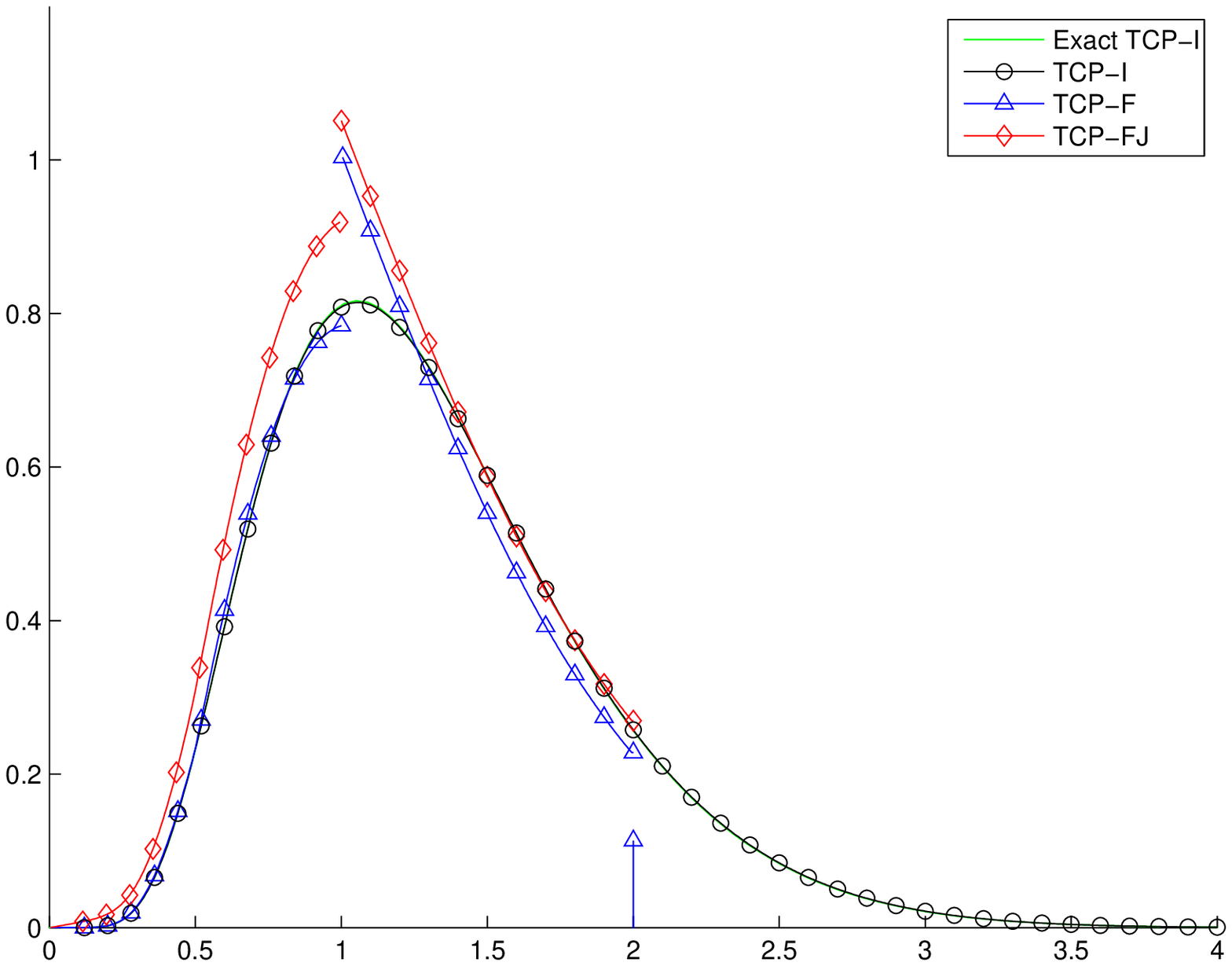}
\quad
\IGS{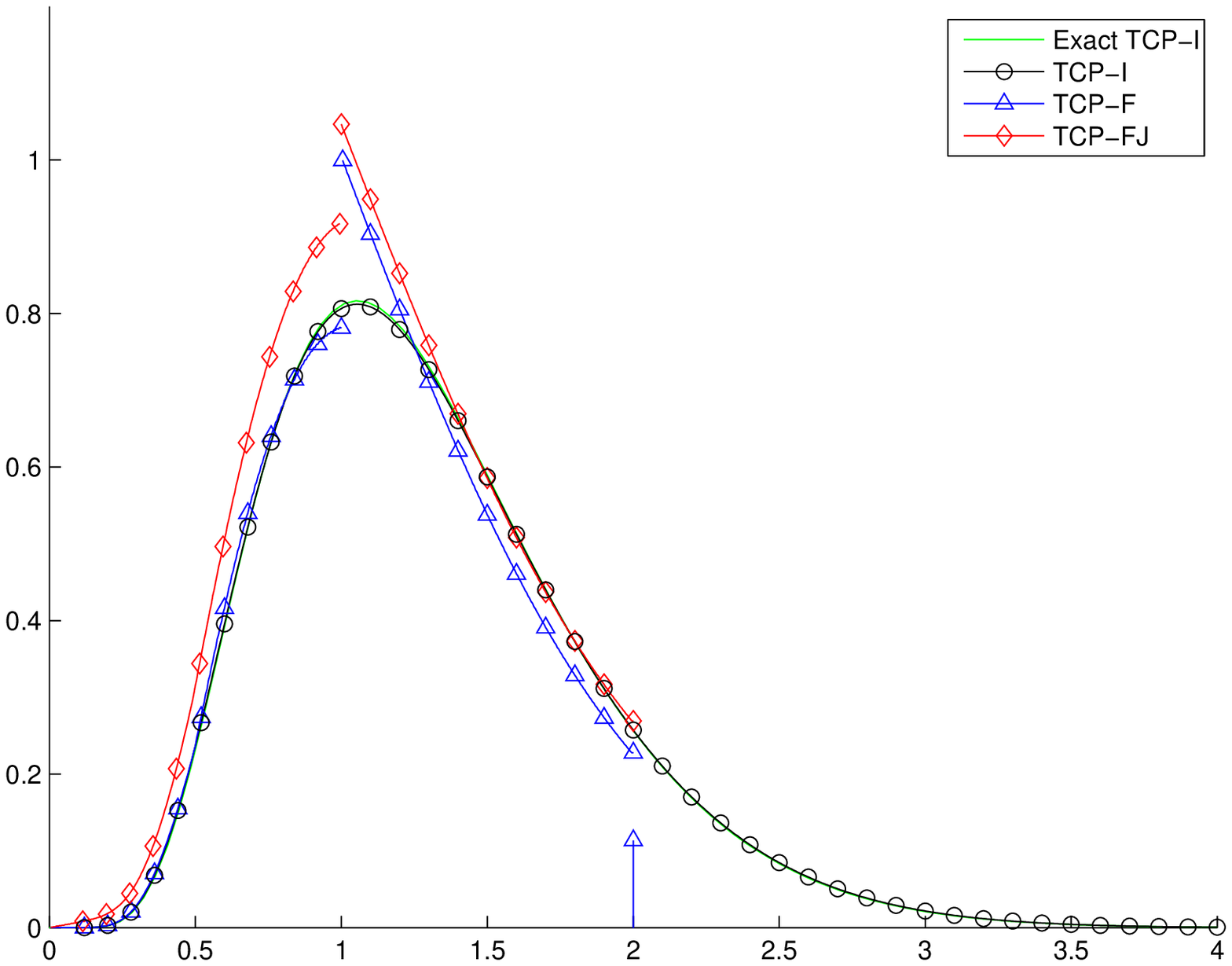}
\caption{Densities for parameters (top-left) $h=\tau=\dt=0.2$, (top-right) $h=\tau=\dt=0.1$, (bottom-left) $h=\tau=\dt=0.01$ and (bottom-right) $h=2\tau=\dt=0.01$.}\label{Fig1}
\end{figure}

In \cite{CMP} the authors get quantitative estimates for the convergence to equilibrium which suggest that an exponential rate could be obtained in many variant of the TCP window size process. So we present the density at time $T=10$ (the equilibrium is reached). For the simulation we have chosen $X=6$ in the TCP-I model, $X=2$ in the TCP-F or TCP-FJ models and $p=0.5$.
We want to illustrate convergence as $h, \tau, \dt\rightarrow 0$.
So we apply the numerical method with the mesh $\cup_{K\leq X/h} [Kh-h,Kh[$ for different parameters $h$, $\tau$ and $\dt=0.01$. Moreover we have drawn in solid line the explicit density of the TCP-I model given in \cite{DGR}. Finally, for the TCP-F model, the asymptotic density has a Dirac mass at point $X$ which is drawn as a vertical arrow whose height is its mass.
The results are illustrated on Figure \ref{Fig1} on a common state space $[0,4]$. 

To conclude we say that the numerical method developed here for Piecewise Deterministic Markov Processes must now be tested on different PDMPs, some of them related to practical cases and others in order to understand their behavior. Indeed the method permits to view the global behavior of the density, whereas the classical approach with particles and Monte-Carlo simulations only describes the trajectories. With our method we can see how the stochastic jumps interact with the flow, which leads sometimes to some unexpected behaviors (see \cite{BLBMZ}). Moreover the method is robust with large values of $\dt$, which permits to compute the asymptotic stationary states of a PDMP.
\medskip

Acknowledgement: The author is thankful for F. Nabet's advice about finite volume simulations.

\end{document}